\documentclass{article}
\usepackage{amsfonts,latexsym}
 
\newcommand{\R}{\mathbb R} 
 
\newcommand{\N}{\mathbb N} 
  
\newcommand{\del}{\partial}
\newcommand{\e}{\varepsilon}

\def\ds{\displaystyle}

\newtheorem{theorem}{Theorem} 
\newtheorem{lemma}{Lemma} 
\newtheorem{proposition}{Proposition} 
 
\newtheorem{definition}{Definition}

\begin{document} 
 
\title{Higher dimensional Scherk's hypersurfaces}
 
\author{Frank Pacard\thanks{Universit\'e Paris 12, 61 Avenue du 
G\'en\'eral de Gaulle, 94010, Cr\'eteil. France.}}

\maketitle 

\begin{abstract}
In $3$-dimensional Euclidean space, Scherk second surfaces are singly 
periodic embedded minimal surfaces with four planar ends. In this paper, 
we obtain a natural generalization of these minimal surfaces in any higher 
dimensional Euclidean space ${\R}^{n+1}$, for $n \geq 3$. More 
precisely, we show that there exist $(n-1)$-periodic embedded minimal 
hypersurfaces with four hyperplanar ends. The moduli space of these 
hypersurfaces forms a $1$-dimensional fibration over the moduli space 
of flat tori in ${\R}^{n-1}$. A partial description of the boundary of this 
moduli space is also given.
\end{abstract}

\section{Introduction}

In $3$-dimensional Euclidean space Scherk second surfaces come in a 
$1$-parameter family $(S_\e)_{\e \in (0,  \frac{\pi}{2})}$ which can be
 described in many different ways. For example it can be described 
{\it via} its Weierstrass representation data  \cite{Die-Hil}, \cite{Kar}
\[
X_\e (\omega) : = \Re \, \int_{\omega_0}^\omega \left( \frac{1}{2}
\left( \frac{1}{g}- g\right), \frac{i}{2}\, \left( \frac{1}{g}+g\right), 1
\right) \, dh_\e ,
\]
where
\[
g (\omega) : = \omega \qquad \mbox{and} \qquad dh_\e := 4 \,\sin \e \, 
(\omega^4 + 1 - 2 \, \cos \e \, \omega^2)^{-1}\, d\omega .
\]
Or even more simply as the zero set of the function
\begin{equation}
F_\e (x_1, x_2, z)  :=  (\cos \e)^2 \, \cosh \left( \frac{x_1}{\cos \e}
\right) - (\sin \e)^2 \, \cosh \left( \frac{z}{\sin \e}\right) - \cos x_2 .
\label{eq:1-1}
\end{equation}
Indeed, it is well known that, the zero set of a function $F$ is a minimal 
surface if and only if $0$ is a regular value of $F$ and
\[
\mbox{div} \left( \frac{\nabla F}{|\nabla F|}\right) =0 ,
\]
on the zero set of $F$. Using this, it is straightforward to check that the 
zero set of $F_\e$ is a minimal surface.

In any of these descriptions, the parameter $\e$ belongs to $(0, \pi/2)$. 
Observe that we do not consider any dilation, translation or rotation of 
a minimal surface, in other words we are only interested in the space 
of surfaces modulo isometries and dilations. Now, we would like to point 
our a few properties of Scherk's second surfaces which will enlighten our 
construction of their higher dimensional analogues. 

\medskip

{\bf (i) - Periodicity~:} Observe that Scherk's second surfaces are singly 
periodic and, in the above description, their common period has been 
normalized to be equal to $(0,2\pi,0)$. Hence, if we define $T^1 : = {\R}/ 
2\pi {\mathbb Z}$, we can consider $S_\e$ to be a minimal surface 
embedded in $ {\R} \times T^1 \times {\R}$. 

\medskip

{\bf (ii) - Asymptotic behavior as $\e$ tends to $0$~:} Another feature 
which will be very important for us is the study the behavior of Scherk's 
second surfaces as the parameter $\e$ tends to $0$ (a similar analysis 
can be performed when the parameter $\e$ tends to $\pi/2$). To this aim, 
we write for all $(x_1,x_2)$ in some fixed compact subset ${\R}^2 - \{0\} 
\times 2\, \pi \,  {\mathbb Z} $ and for all $\e$ small enough
\[
z = \pm \sin \e \, \mbox {acosh} \left(  (\tan \e)^{-2} \, \cosh \left( 
\frac{x_1}{\cos \e}\right) - (\sin  \e )^{-2}\, \cos x_2 \right).
\]
Using this, we readily see that, away from the set $\{0\} \times 2\, \pi 
{\mathbb Z} $, the one parameter family of surfaces $S_\e$ converges 
to the union of two horizontal planes, as $\e$ tends to $0$. In other words,
 the sequence of surfaces $S_\e$ converges, away from the origin, 
to two copies of ${\R} \times T^1\times \{0\}$ in ${\R}\times T^1 \times 
{\R}$, as the parameter $\e$ tends to $0$.

As already mentioned, a similar analysis can  be carried out as the 
parameter $\e$ tends to $\pi/2$ and, this time, we find that the sequence 
of surfaces $S_\e$ converges, away from the origin, to two copies of 
$\{0\} \times T^1 \times {\R}$ in ${\R}\times T^1 \times {\R}$.

\medskip

{\bf (iii) - Blow down analysis~:} For each fixed $\e \in (0, \pi/2)$, the 
surface $S_\e$ has four planar ends which are asymptotic to 
\[
V_\e^\pm : = \left\{ (x_1,x_2,z) \in{\R} \times  T^1 \times {\R} \, : \, z = 
\pm \, \left( \tan \e \, |x_1| - 2 \, \sin \e \, \log \tan \e \right) \right\}.
\]
More precisely, away from a compact set in ${\R}\times T^1 \times 
{\R}$, the surface $S_\e$ is a normal graph over $V_\e^\pm$ for some 
function which is exponentially decaying as $x_1$ tends to $\pm \infty$. 
Another way to understand this would be to say that, the sequence of 
surfaces $\lambda \, S_\e$ converges, as $\lambda$ tends to $0$ to 
$W_\e^+ \cup W^-_\e$, where 
\[
W_\e^\pm : =\left\{ (x_1,x_2,z) \in {\R}\times T^1 \times {\R} \, : \, z 
= \pm \, \tan \e \, |x_1|  \right\}.
\]

\medskip

{\bf (iv) - Blow up analysis~:} Instead of blowing down the surfaces 
$S_\e$ as we have done in (iii), we can blow up the surfaces $S_\e$ 
by considering the sequence of scaled surfaces $\e^{-1} \, S_\e$. As 
$\e$ tends to $0$ this sequence converges on compact to a vertical 
catenoid. To see this, just define the new set of coordinates
\[
(\tilde x_1,\tilde x_2, \tilde z) : = \frac{1}{2 \sin \e} (x_1,x_2,z),
\]
and, in  (\ref{eq:1-1}), we expend both $\cos x_2$ and $\cosh (x_1 / \cos 
\e)$, in terms of powers of $\e$. We find with little work
\[
(\cos \e)^2 \, \left( 1+ 2 \, (\tan \e)^2 \, \tilde x_1^2\right) - (\sin \e)^2 \, 
\cosh (2 \, \tilde z)  = 1 - 2 \, (\sin \e)^2 \, \tilde x_2^2 + {\cal O}( \e^4) .
\]
Hence 
\[
\tilde x_1^2+  \tilde x_2^2 =  \cosh^2  \tilde z + {\cal O}(\e^2) .
\]
Clearly, as $\e$ tends to $0$, this converges, uniformly on compact sets, 
to an implicit parameterization of a  vertical catenoid.

\medskip

To complete this brief description, let us mention that Scherk's second 
surfaces have recently been used as one of the building blocks of some 
desingularization procedure, to produce new embedded minimal surfaces 
in $3$-dimensional Euclidean space. We refer to the work of M. Traizet 
\cite{Tra-1} and also to the recent work of N. Kapouleas \cite{Kap-1}, 
\cite{Kap-2} for further details.

\medskip

In order to state our result properly, we need to introduce two ingredients 
which will be fundamental in our analysis.  First observe that, in higher 
dimensions, there is a natural generalization of the catenoid in 
Euclidean $3$-space. This hypersurface, which we will call the unit 
$n$-catenoid, is a hypersurface of revolution with two hyperplanar ends. 
It can be parameterized by 
\[
{\R} \times S^{n-1} \ni (s, \theta)  \longrightarrow (\varphi (s) \, \theta, 
\psi (s))\in  {\R}^{n+1},
\]
where the function $\varphi$ is defined by the identity $\varphi^{n-1} 
(s) = \cosh ((n-1) s))$ and where the function $\psi$ is given by
\[
\psi (s) : = \int_0^s  \varphi^{2-n}(t) \, dt .
\]
Using this $n$-catenoid, S. Fakhi and the author have produced examples 
of complete immersed minimal hypersurfaces of ${\R}^{n+1}$ which 
have $k\geq 2$ hyperplanar ends \cite{Fak-Pac}. These hypersurfaces 
have the topology of a sphere with $k$ punctures and they all have finite 
total curvature, they generalize the well known $k$-noids in $3$-dimensional 
Euclidean space \cite{Hof-Mee}.
 
\medskip

Another ingredient in our analysis is the moduli space of flat tori in 
${\R}^{m}$, for $m \geq 1$. We recall a few well known facts about this 
moduli space and refer to \cite{Wol} for further details. Any flat torus in 
${\R}^m$ can be identified with ${\R}^m / A {\mathbb Z}^m$ where 
$A\in GL(m,{\R})$. The volume of the $m$-dimensional torus $T^m : = 
{\R}^m / A {\mathbb Z}^m$ is then given by 
\[
\mbox{vol} \, (T^m ) =  | \mbox{det} \, A| .
\]
It is a simple exercise to check that two tori ${\R}^m / A {\mathbb Z}^m$ 
and ${\R}^m / B {\mathbb Z}^m$  are isometric if and only if there exist 
$M \in O(m, {\mathbb R})$ and $N \in GL (m , {\mathbb Z})$ such that 
$A = M \, B \, N$. The moduli space of flat tori ${\cal T}^m$ is defined 
to be the space of flat tori ${\R}^m / A {\mathbb Z}^m$, normalized by 
asking that 
\[
\mbox{vol} \, (T^{m}) = \mbox{vol}\,  (S^{m}),
\]
modulo isometries. For later use, it will be convenient to identify any torus 
$T^m \in {\cal T}^m$ with a subset of ${\R}^m$. To this aim, if
\[
T^m ={\R}^{m}/ A {\mathbb Z}^m ,
\]
for some $A \in GL(n , {\R})$, we identify $T^m$ with the image of 
$[-\frac{1}{2}, \frac{1}{2}]^m$ by $A$. In particular, we will talk about 
the origin $0 \in T^m$, simply referring to the origin in $A \, [-\frac{1}{2}, 
\frac{1}{2}]^m \subset {\R}^m$.  We will also consider $B^n_\rho \subset 
{\R}^{n-m}\times T^m$ as the $n$-dimensional ball of radius $\rho$ in 
${\R}^{n-m}\times A \, [-\frac{1}{2}, \frac{1}{2}]^m$. And so on. 
Also observe that, granted this identification, $T^m$ is invariant under the action 
of the following subgroup of $O(m , {\mathbb R})$
\[
{\mathfrak D}_m : =  \left\{ D : = \mbox{diag} \, (\eta_1, \ldots , \eta_m) 
\quad | \quad \eta_i = \pm 1 \right \} .
\]

\medskip

In this paper, we pursue the quest of higher dimensional generalizations 
of classical minimal surfaces which we have initiated in \cite{Fak-Pac}. 
More precisely, we obtain a natural generalization of Scherk's 
second surfaces in higher dimensional Euclidean spaces. Recall that one 
can view the moduli space of Scherk's surfaces as a $1$-dimensional 
fibration over the moduli space of flat tori in ${\R}$. We will show that,
 in ${\R}^{n+1}$, for $n \geq 3$, there exists a finite dimensional family 
of embedded minimal hypersurfaces satisfying properties which are similar 
to (i)-(iv). This family, which turn out to be a $1$-dimensional fibration 
over the moduli space of flat tori in ${\R}^{n-1}$, yields a partial description 
of the moduli space of what might be called "higher dimensional Scherk's 
hypersurfaces". More precisely, we obtain a description of the boundary 
of this moduli space, this boundary turns out to be modeled over the moduli 
space of tori in ${\R}^m$ for any $1\leq m \leq n-1$.  

\medskip

Our main result can be stated as follows~:
\begin{theorem}
Assume that $n \geq 3$ and $1 \leq m \leq n-1$ are fixed. Let $T^{m} \in 
{\cal T}^{m}$ be any flat torus of ${\R}^{m}$.  Then, there exist $\e_0 >0$ 
and $(S_\e)_{\e \in (0, \e_0)}$ a one parameter family of minimal 
hypersurfaces of ${\R}^{n-m}\times T^{m}\times {\R}$ such that~:
\begin{enumerate}

\item[(i)]  For all $\e \in (0, \e_0)$, the hypersurface $S_\e$ is embedded in 
${\R}^{n-m}\times T^{m}\times {\R}$ and is invariant under the action of  
$O(n-m, {\R}) \otimes {\mathfrak D}_m \otimes  \{\pm I_1\}\subset O(n+1, 
{\R})$.

\item[(ii)] As $\e$ tends to $0$, the sequence of hypersurfaces $(S_\e)_\e$ 
converges to the union of two copies of $ {\R}^{n-m} \times T^{m}\times 
\{0\}$, away from the origin.

\item[(iii)] For all $\e \in (0, \e_0)$, there exists $c_\e >0$ and 
$d_\e >0$ such that the 
hypersurface $S_\e$ has four ends which are asymptotic to
\[
V_\e^\pm : = \{(x_1, x_2,z) \in {\R}^{n-m} \times T^{m} \times  {\R} 
\quad : \quad z = \pm  \, ( c_\e \, \zeta_m (x_1) + d_\e ) \},
\]
where $\zeta_{n-1}(y) := |y|$, $\zeta_{n-2}(y) := \log |y|$
 and $\zeta_m (y) := 0$, when $m \leq n-3$. In particular, this means that, 
up to a translation along the $z$-axis, the hypersurface $S_\e$ is a normal 
graph over $V_\e^\pm$ for some function which is exponentially 
decaying in $|x_1|$ function when $m=n-1$ and for some function which 
is polynomially decaying in $|x_1|$ when $m\leq n-2$. Furthermore, when 
$m=n-1$, we have 
\begin{equation}
\lim_{\e \rightarrow 0} \frac{c_\e}{\e^{n-1}} = \frac{1}{2} .
\label{eq:slope}
\end{equation}

\item[(iv)] As $\e$ tends to $0$, the sequence of rescaled hypersurfaces 
$( \e^{-1} \, S_\e)_\e$ converges, uniformly on compact sets, to a vertical 
unit $n$-catenoid. 
\end{enumerate}
\end{theorem}

This result, when $m = n-1$, yields minimal hypersurfaces which constitute 
the natural generalization of Scherk's second surfaces in higher dimensional 
Euclidean spaces. More precisely, when $m =n -1$, the above result provides 
a description of part of ${\cal S}_{n}$, the moduli space of $n$-dimensional 
Scherk's hypersurfaces in ${\R}^{n+1}$. As this result shows, this moduli 
space is locally a $1$-dimensional fibration over the moduli space of flat tori 
in ${\R}^{n-1}$. Though we have not been able to prove it, we expect this 
fibration to extend, as it does when $n=1$, to all $c_\e \in (0, \pi/2)$. 

The above result, when $m \leq n-2$, yields hypersurfaces which have to be 
understood as belonging to the boundary of the moduli space ${\cal S}_{n}$, 
in the same way that any product $ {\R}^{n-m -1} \times T^m$, for $m \leq 
n-2$ corresponds to a point in the compactification of the moduli space of flat 
tori in ${\R}^{n-1}$.  We expect that the moduli space ${\cal S}_n$ can be 
compactified and that the family of hypersurfaces described in the above result 
constitute a collar neighborhood of the boundary of $\overline {\cal S}_{n}$. 
In other words, Theorem~1 should provide a local description of ${\cal S}_{n}$,
 near its boundary.

\medskip

To  conclude, let us briefly describe the strategy of the proof of the result. It 
should be clear from (ii) and (iv) that, for small $\e$, Scherk's second surfaces 
can be understood as a desingularization of two copies of ${\R}\times T^1 
\times \{0\}$ in ${\R}\times T^1\times {\R}$. Keeping this observation in mind, 
our strategy will be to show that a similar desingularization is possible for two 
copies of ${\R}^{n-m} \times  T^{m}\times \{0\}$ in ${\R}^{n-m} \times  
T^{m}\times {\R}$. The proof of this result is very much in the spirit of 
\cite{Fak-Pac}, \cite{Maz-Pac} or \cite{Maz-Pac-Pol}, however some aspects 
are simpler in the present paper thanks to the special geometry of our problem.

\medskip

Our work has been strongly influenced by the recent work of M. Traizet 
\cite{Tra-2} and the work of N. Kapouleas \cite{Kap-1}, \cite{Kap-2} in their 
construction of minimal embedded surfaces in ${\R}^3$. Indeed, on the one
 hand, N. Kapouleas has used Scherk's second surfaces to desingularize finitely 
many catenoids or planes having a common axis of revolution and he has produced 
embedded minimal surfaces with finitely many ends and very high genus. On 
the other hand, M. Traizet has used finitely many catenoids to desingularized  
parallel planes and produced minimal surfaces with finitely many ends and any 
genus (larger than $2$).  There is a formal link between these two constructions 
since, in some vague sense,  the surfaces constructed by N. Kapouleas on the 
one hand and the surfaces constructed by M. Traizet, for a genus large enough, 
on the other hand, should belong to the same moduli space. It was therefore 
tempting to try to produce Scherk's second surfaces using some desingularization 
procedure.

\medskip

\noindent
{\bf Acknowledgment~:} This paper was written while author was visiting the 
Mathematical Sciences Research Institute in Berkeley. He would like to take 
this opportunity to thank the MSRI for their support and hospitality. The author 
would also like to thank  D. Hoffman, R. Mazzeo and M. Weber, for stimulating 
discussions.

\section{Definitions and notations}

In this brief section we record some notations and definitions which will be used 
throughout the paper.  

\medskip

{\bf Eigenfunctions of $\Delta_{T^{m}}$~:} Given $m \geq 1$ and $T^{m}\in 
{\cal T}^{m}$, we will denote by $E_i$, $i\in {\mathbb N}$ the eigenfunctions 
of the Laplacian on $T^{m}$ with corresponding eigenvalues $\mu_j$, that is 
$\Delta_{T^m} E_i = - \mu_i \, E_i$, with $\mu_i \leq \mu_{i+1}$. We will 
assume that these eigenfunctions are counted with multiplicity and are normalized 
so that
\[
\int_{T^m} E_i^2 \, dx =1.
\]
Though the spectral data of $\Delta_{T^m}$ do depend on $T^m$, we will not 
write this dependence in the notation.

\medskip

{\bf Functions on $T^{m}$ which are invariant under the action of some group~:}
 We will only be interested in function on $T^m$ and eigenfunctions of 
$\Delta_{T^{m}}$ which have some special symmetry. Namely the set of 
functions and eigenfunctions which are invariant under the action of the following 
subgroup of $O(m , {\mathbb R})$
\[
{\mathfrak D}(m) :=  \left\{ D : = \mbox{diag} \, (\eta_1, \ldots, , \eta_m) \quad 
| \quad \eta_\ell = \pm 1 \right \} .
\]
We define ${\mathfrak  I}(m) \subset {\mathbb N}$ to be the set of indices $i$ 
corresponding to eigenfunctions $E_i$ which are invariant under the action of $
{\mathfrak D}(m)$, that is
\begin{equation}
{\mathfrak  I} (m) : = \left\{i \geq 0 \quad | \quad  E_i = E_i \circ D, \quad  
\mbox{for all} \quad  D \in {\mathfrak D}(m) \right\}.
\label{eq:2-1}
\end{equation}

\medskip

{\bf Eigenfunctions of $\Delta_{S^{n-1}}$~:} For all $n \geq 2$, we will denote 
by $e_j$, $j \in {\mathbb N}$, the eigenfunctions of the Laplacian on $S^{n-1}$ 
with corresponding eigenvalues $\lambda_j$, that is $\Delta_{S^{n- 1}} e_j 
= - \lambda_j \, e_j$, with $\lambda_j \leq \lambda_{j+1}$. We will assume 
that these eigenfunctions are counted with multiplicity and are normalized so that
\[
\int_{S^{n-1}} e_j^2 \, d\theta =1.
\]

\medskip

{\bf Functions on ${\R}^n$ or on $S^n$ which are invariant under the action 
of some group~:} Given $1\leq m \leq n-1$, we can decompose ${\R}^{n} =  
{\R}^{n-m} \times {\R}^{m}$. We will only be interested in function on ${\R}^n$ 
and eigenfunctions of $\Delta_{S^{n-1}}$ which have some special symmetry. 
Namely functions which are invariant under the action of the following subgroup 
of $O(n, {\R})$
\[
{\mathfrak H}(n,m) :=  O(n -m , {\mathbb R})  \otimes {\mathfrak D}(m) .
\]
It will be convenient to define ${\mathfrak J}(n,m)$ to be the set of indices 
$j \in {\mathbb N}$ corresponding to eigenfunctions $e_j$ which are invariant 
under the action of ${\mathfrak H}(n,m)$, that is
\[
{\mathfrak J}(n,m) : = \left\{  j \geq 0 \quad | \quad  e_j = e_j \circ R,  \quad  
\mbox{for all} \quad  R \in {\mathfrak H}(n,m)  \right\}.
\]
It will be important to observe that $1, 2, \ldots, n$ do not belong to 
${\mathfrak J}(n,m)$ since the eigenfunctions corresponding to the eigenvalues 
$\lambda_1= \ldots = \lambda_n$ are not invariant under the action of 
$-I_n  \in {\mathfrak H}(n,m)$. 

\medskip

For all $k \in {\N}$ and all $\alpha \in (0,1)$, we define ${\cal C}^{k, \alpha}
(S^{n-1}, {\mathfrak H}(n,m))$ to be the subset of functions of ${\cal C}^{k, 
\alpha} (S^{n-1})$ whose eigenfunction decomposition only involves indices 
belonging to ${\mathfrak J}(n,m)$. In other words,  $g \in {\cal C}^{k, \alpha} 
(S^{n-1}, {\mathfrak H}(n,m))$ if and only if $g \in {\cal C}^{k, \alpha}
 (S^{n-1})$ and 
\[
g =  g_0 + \sum_{j \in {\mathfrak J}} g_j \, e_j .
\]
Observe that, by definition, any function of ${\cal C}^{k, \alpha} (S^{n-1}, 
{\mathfrak H}(n,m))$ is orthogonal to $e_1, \ldots, e_{n-1}$ in the $L^2$ sense, 
on $S^{n-1}$. 

\medskip

{\bf Notations~:} Given $1\leq m \leq n-1$, we will adopt the following notation~: 
\[
x \quad  \mbox{or} \quad (x_1, x_2) \in {\R}^{m}\times {\R}^{m} \sim 
{\R}^{n-m} ,
\]
will denote a point in ${\R}^n$ and  
\[
(x,z) \in {\R}^{n}\times {\R} \sim {\R}^{n+1} ,
\]
will denote a point in ${\R}^{n+1}$. Finally, $\theta$ will denote a  point in 
$S^{n-1}$.

\section{Minimal hypersurfaces close to a truncated $n$-catenoid}

This section is mainly adapted from \cite{Fak-Pac}, we first recall some of 
the technical results of \cite{Fak-Pac} which are needed in this paper and 
adapt them to our situation. 

\subsection{The $n$-catenoid}

Assume that $n \geq 3$ is fixed. We recall some well known fact concerning 
the unit $n$-catenoid $C_1$ which is a minimal hypersurface of revolution in 
${\R}^{n+1}$, further details are available in \cite{Fak-Pac}. By definition, 
 $C_1$ is the minimal hypersurface of revolution parameterized by 
\begin{equation}
X_0 : (s, \theta) \in {\R} \times S^{n-1}  \longrightarrow (\varphi (s) \, \theta, 
\psi (s))\in  {\R}^{n+1},
\label{eq:3-6}
\end{equation}
where  $\varphi$ to be the unique, smooth, non constant solution of 
\[
(\del_s \varphi)^2 + \varphi^{4-2n} =\varphi^2 \qquad \mbox{with}\qquad 
\varphi(0)=1 ,
\]
and where the function $\psi$ is the unique solution of
\[
\del_s \psi = \varphi^{2-n} \qquad \mbox{with}\qquad \psi(0)=0.
\]
As already mentioned in the introduction, it might be interesting to observe that 
$\varphi$ is explicitely given by the identity
\[
\varphi^{n-1} (s) = \cosh ((n-1) s)) .
\]
Using this, it is easy to check that the function $\psi$ converges at $\pm \infty$. 
We set 
\[
c_\infty : = \lim_{s\rightarrow +\infty} \psi .
\]
The fact that $\psi$ converges at both $\pm \infty$ implies that the hypersurface 
$C_1$ has two hyperplanar ends and is in fact contained between the two 
asymptotic hyperplanes defined  by $z =  \pm c_\infty$. In addition,  the upper 
end (resp. lower end) of the unit $n$-catenoid can be parameterized as a graph 
over the $z = 0$ hyperplane for some function $u$ (resp $-u$). It is easy to 
check that the function $u$ has the following expansion as $r : =  |x|$ tends to 
$\infty$
\begin{equation}
u =  c_\infty  -  \frac{1}{n-2} \, r^{2-n} + {\cal O}( r^{4-3n}).
\label{eq:dldl}
\end{equation}

\subsection{The mean curvature operator}

Let us assume that the orientation of $C_1$ is chosen so that the unit normal 
vector field is given by 
\begin{equation}
N_0 : = \frac{1}{\varphi}\left( \del_s \psi \, \theta, - \del_s \varphi \right) .
\label{eq:3-7}
\end{equation}
All surfaces close enough to $C_1$ can be parameterized (at least locally) as 
normal graphs over $C_1$, namely as the image of
\[
X_w  : = X_0  + w \, \varphi^{\frac{2-n}{2}} \, N_0,
\]
for some small function $w$. The following technical result is borrowed from  
\cite{Fak-Pac}. It just states that the mean curvature of the hypersurface 
parameterized by $X_w$ has some nice expansion in terms of $w$. Observe 
that, in order to define $X_w$, we have used $w \, \varphi^{\frac{2-n}{2}} \, 
N_0$ instead of the usual $w \,  N_0$, there is no loss of generality in doing 
so and this choice will simplify the notations in the forthcoming result which 
describes the structure of the nonlinear partial differential equation $w$ has 
to satisfy in order for the hypersurface parameterized by $X_w$ to be minimal.
\begin{proposition} \cite{Fak-Pac}
The hypersurface parameterized by $X_w$ is minimal if and only if the 
function $w$ is a solution of the nonlinear elliptic partial differential equation
\begin{equation}
{\cal L}\,  w = \varphi^{\frac{2-n}{2}}\, Q_2 \left(\varphi^{-\frac{n}{2}}\, w 
\right)+ \varphi^{\frac{n}{2}} \, Q_3 \left( \varphi^{-\frac{n}{2}}\, w \right),
\label{eq:3-1}
\end{equation}
where 
\[
{\cal L} : = \del_{s}^2+ \Delta_{S^{n-1}} -\left( \frac{n-2}{2}\right)^2 + 
\frac{n \, (3n-2)}{4}\, \varphi^{2-2n} ,
\]
where $\xi \rightarrow Q_2(\xi)$ is a nonlinear second order differential operator 
which is homogeneous of degree $2$ and where $\xi \rightarrow Q_3(\xi)$ is a 
nonlinear second order differential operator which satisfies
\[
Q_3 (0)=0, \quad D_\xi Q_3(0)=0\quad \mbox{and}\quad D^2_\xi Q_3(0)=0 .
\]
Furthermore, the coefficients of $Q_2$ on the one hand and the coefficients in 
the Taylor expansion of $Q_3$ with respect to the $\xi$, computed at any $\xi$ 
in some fixed  neighborhood of $0$ in ${\cal C}^{2, \alpha}({\R}\times 
S^{n-1})$ on the other hand are bounded functions of $s$ and so are the 
derivatives of any order of these functions. 
\label{pr:1}
\end{proposition}
The operator ${\cal L}$ is clearly equivariant with respect to any action of the form
\[
 {\R}\times S^{n-1} \ni (s, \theta)  \longrightarrow (-s, R\, \theta)  \in {\R}\times S^{n-1}, 
\]
when $R \in {\mathfrak H}(n,m)$. Since in addition the mean curvature is invariant 
by isometries, we conclude that the nonlinear operator which appears on the right 
hand side of (\ref{eq:3-1}) also enjoys this equivariance property.

\medskip

It might be useful to rephrase the properties of the nonlinear operators $Q_2$ and 
$Q_3$ into a slightly weaker form. It follows from the properties of $Q_2$ and 
$Q_3$ that there exist constants $c, c_0 >0$ such that, for all $s \in {\R}$ and all 
$\xi_1, \xi_2 \in {\cal C}^{2, \alpha}([s,s+1]\times S^{n-1})$, we have 
\begin{equation}
\| Q_2 (\xi_1)- Q_2 (\xi_2)\|_{{\cal C}^{0, \alpha}} \leq c \, \left( \sup_{i=1,2}
\|\xi_i \|_{{\cal C}^{2, \alpha}}\right) \,  \|\xi_2 -\xi_1 \|_{{\cal C}^{2, \alpha}}
\label{eq:useful1}
\end{equation}
and, provided $\|\xi_i \|_{{\cal C}^{2, \alpha}}\leq c_0$, we also have 
\begin{equation}
\| Q_3 (\xi_1)- Q_3 (\xi_2)\|_{{\cal C}^{0, \alpha}} \leq c \, \left( \sup_{i=1,2}
\|\xi_i \|_{{\cal C}^{2, \alpha}}\right)^2  \, \|\xi_2 -\xi_1 \|_{{\cal C}^{2, \alpha}},
\label{eq:useful2}
\end{equation}
where all norms are understood on $[s,s+1]\times S^{n-1}$. Since $Q_2$ is 
homogeneous of degree $2$ no assumptions on $\xi_i$ are required in order to get 
the estimate involving $Q_2$, however they are required for the estimates involving
 $Q_3$.

\medskip

Let us warn the reader that the operator ${\cal L}$ which appears in this result is 
not the Jacobi operator which is defined to be the linearized mean curvature operator
 when nearby hypersurfaces are normal graphs over the $n$-catenoid, that is when 
they are parameterized by 
\[
\tilde X_w  : = X_0  + w \, N_0,
\]
but ${\cal L}$ is conjugate to the Jacobi operator.

\subsection{Linear analysis}

Projecting the operator ${\cal L}$ over the eigenspaces spanned by $e_j$, we are left 
with the study of the sequence of operators
\[
L_j : = \del_{s}^2 -\lambda_j - \left(\frac{n-2}{2}\right)^2 +\frac{n (3n-2)}{4}\, 
\varphi^{2-2n}, \qquad j \in {\mathbb N}.
\]
The indicial roots of ${\cal L}$ at both $+\infty$ or $-\infty$ are given by $\pm 
\gamma_j$ where
\begin{equation}
\gamma_j : = \sqrt{\left(\frac{n-2}{2}\right)^2+\lambda_j}.
\label{eq:3-4}
\end{equation}
Let us recall that these indicial roots appear in the study of the asymptotic behavior 
at $\pm \infty$ of the solutions of the homogeneous problem $L_j w =0$. More 
precisely, for each $j \in {\mathbb N}$, one can find $w^{\pm}_j$, two independent 
solutions of $L_j w =0$ such that $w_j^+ (s) \sim e^{\gamma_j s}$ and $w_j^-( s)
\sim e^{- \gamma_j s}$ at $+ \infty$. Observe that the functions $s \rightarrow 
w_j^\pm (-s)$ are solutions of $L_j w =0$ such that $w_j^+ (s) \sim e^{-\gamma_j s}$  
and $w_j^- (s) \sim e^{\gamma_j s}$ at $- \infty$. These indicial roots will play a 
crucial r\^ole in the study of the mapping properties of ${\cal L}$.

\medskip

\noindent
To keep the notations short, we define the second order elliptic operator
\[
\Delta_0 : = \del_s^2 + \Delta_{S^{n-1}}- \ds\left(\frac{n-2}{2}\right)^2 ,
\]
which acts on functions defined on $\R\times S^{n-1}$. In particular
\[
{\cal L} =  \Delta_0  +\frac{n (3n-2)}{4}\, \varphi^{2-2n} .
\]
The indicial roots of $\Delta_0$ at both $+\infty$ or $-\infty$ are also given by $\pm 
\gamma_j$.

\medskip

It is straightforward to check that $\Delta_0$ satisfies the maximum principle and 
also that the operator ${\cal L}$ does not satisfy the maximum principle because 
of the presence of the extra potential. Indeed, one can check that the functions
\begin{equation}
\Psi^{0,-}: = \del_s ( \varphi^{\frac{n-2}{2}} ) , \qquad \qquad  
\Psi^{0,+}: = \varphi^{\frac{n-4}{2}} \, (\varphi \, \del_s \psi- \psi \, \del_s \varphi ),
\label{eq:jf0}
\end{equation}
and, for $j=1, \ldots, n$,  the functions
\begin{equation}
\Psi^{j, -}: = \varphi^{\frac{n-4}{2}} \, (\varphi \, \del_s \varphi + \psi \, \del_s \psi ) \, e_j,
\qquad \qquad  
\Psi^{j, +}: = \varphi^{-\frac{n}{2}} \, e_j ,
\label{eq:jf}
\end{equation}
are  Jacobi fields, i.e. are solutions of the homogeneous problem ${\cal L} w =0$, 
and that the $\Psi^{j, +}$ are bounded. Nevertheless, the following result, borrowed 
from \cite{Fak-Pac},  asserts that, the operator ${\cal L}$ still satisfies the 
maximum principle if it is restricted to the higher eigenspaces of the cross-sectional 
Laplacian $\Delta_{S^{n-1}}$~:
\begin{proposition} 
Assume that $\delta < \frac{n+2}{2}$ is fixed and that $w$ is a solution of 
\[
{\cal L} w = 0 ,
\]
 which is bounded by $\varphi^{\delta}$ on $(s_1, s_2)\times S^{n-1}$ and which 
satisfies $w=0$ on $\{s_i\}\times S^{n-1}$, if any of the $s_i$ is finite. Further 
assume that, for each $s \in (s_1, s_2)$, the function  $w(s, \cdot)$ is orthogonal 
to $e_0, \ldots, e_n$ in the $L^2$ sense on $S^{n-1}$. Then $w\equiv 0$.
\label{pr:2}
\end{proposition}

\medskip

In view of the previous result, it is natural to consider the operator ${\cal L}$ acting 
on functions bounded by a constant times a power of the function $\varphi$. As in 
\cite{Maz-Pac} and \cite{Fak-Pac}, we define a family of weighted H\"older 
spaces by~:
\begin{definition}
For all $\delta \in {\R}$, the space ${\cal C}^{k, \alpha}_{\delta}({\R} \times 
S^{n-1})$ is defined to be the space of functions $w \in {\cal C}^{k, \alpha}_{loc} 
( {\R} \times S^{n-1})$  for which the following norm is finite
\[
\|w\|_{{\cal C}^{k, \alpha}_{\delta}}: = \sup_{s \in {\R}} \,  \varphi^{-\delta} \,  
|w |_{{\cal C}^{k, \alpha}([s, s+1]\times S^{n-1})}.
\]
Here $| \, \, \, |_{{\cal C}^{k, \alpha} ( [s, s+1]\times S^{n-1})}$ denotes the 
H\"older norm in $[s, s+1]\times S^{n-1}$.
\label{de:00}
\end{definition}
Moreover, for any $S >0$, the space ${\cal C}^{2, \alpha}_\delta ([-S,S] \times 
S^{n-1})$ is defined to be the space of restriction of functions of  ${\cal C}^{2, 
\alpha}_\delta ({\R}\times S^{n-1})$ to $[-S,S] \times S^{n-1}$. This space is 
naturally endowed with the induced norm.

\medskip

Though this will not be necessary for the remaining of the analysis, we quote here 
some well known properties of the operator
\[
{\cal L} : {\cal C}^{k, \alpha}_{\delta}({\R} \times S^{n-1}) \longrightarrow 
{\cal C}^{k, \alpha}_{\delta}({\R} \times S^{n-1}) .
\]
To keep track of the weighted space over which ${\cal L}$ is 
defined, we will denote the above operator by ${\cal L}_\delta$. The most important 
fact is that the mapping properties of ${\cal L}_\delta$ crucially depend on the 
choice of the weight parameter $\delta$. Indeed, it follows from general arguments 
that ${\cal L}_\delta$ has close range and is even Fredholm if and only if the weight 
$\delta$ is not equal to any of the indicial roots $\pm \gamma_j$, $j \in 
{\mathbb N}$, (a fact which, given the special structure of our operator, can be 
easily proven be separation of variables). The fact that the functions given in 
(\ref{eq:jf}) are Jacobi fields shows that ${\cal L}_\delta$ is not injective when 
$\delta > - \frac{n}{2}$ and it can be proven, with the help of Proposition~\ref{pr:2}, 
that ${\cal L}_\delta$ is injective if $\delta < - \frac{n}{2}$. This later fact in turn 
implies that ${\cal L}_\delta$ is surjective if $\delta > \frac{n}{2}$ is not equal to 
any $\gamma_j$, $j \geq 0$ (this uses the fact that the operator ${\cal L}_\delta$ 
and ${\cal L}_{-\delta}$ are, in some sense, dual).

\medskip

As already mentioned in \S 2, we will only be interested in functions which are 
invariant under the action of some group.  This is the reason why we introduce the~:
\begin{definition}
For all $k \in {\mathbb N}$, $\alpha \in (0,1)$ and $\delta \in {\R}$, the space 
${\cal C}^{k, \alpha}_{\delta }({\R} \times 
S^{n-1}, {\mathfrak H}(n , m))$ is defined to be the space of functions $w \in 
{\cal C}^{k, \alpha}_\delta ( {\R} \times S^{n-1})$ which satisfy
\[
\forall \, (s, \theta) \in {\R}\times S^{n-1}, \qquad \qquad w(s, \theta ) = 
w(- s, \theta),
\]
and also 
\[
\forall \, (s, \theta) \in {\R}\times S^{n-1}, \qquad \qquad w(s,  R \theta ) = 
w(s, \theta),
\]
for all $R\in {\mathfrak H}(n,m)$. This space is endowed with the induced norm.
\label{de:1}
\end{definition}

\medskip

Observe that, any function $w \in {\cal C}^{k, \alpha}_\delta ({\R}\times S^{n-1}, 
{\mathfrak H}(n,m))$ can be decomposed as 
\[
w (s, \theta) = \sum_{j \in {\mathfrak J}} w_j (s) \, e_j (\theta),
\]
where, for all $j$,  all functions $s \longrightarrow w_j(s)$ are even.

\medskip

Observe that the Jacobi fields defined in (\ref{eq:jf}) are not invariant with respect 
to the action of ${\mathfrak H}(n,m)$, hence one can show that
\[
{\cal L} : {\cal C}^{k, \alpha}_{\delta}({\R} \times S^{n-1}, {\mathfrak H}(n,m) ) 
\longrightarrow {\cal C}^{k, \alpha}_{\delta}({\R} \times S^{n-1}, {\mathfrak H}
(n,m)) 
\]
is injective for all $\delta < \frac{n-2}{2}$ and surjective for all $\delta > 
\frac{2-n}{2}$, which is not equal to any $\gamma_j$, for $j\geq 0$. We will not 
need such general statements but, since we will be working with functions defined 
on $[-S, S]\times S^{n-1}$.

\medskip

Among the Jacobi fields defined in (\ref{eq:jf0}) and (\ref{eq:jf}), 
\[
\Psi^{+,0} = \varphi^{\frac{n-4}{2}} \, (\varphi \, \del_s \psi - \psi \, \del_s \varphi).
\]
is the only one which is invariant with respect to the action of ${\mathfrak H}(n,m)$ 
and which is an even function of $s$. It is easy to see that this Jacobi field vanishes 
for finitely many values of $s$. Let us define $s_0 >0$ to be the largest zero of  the 
function $\Psi^{+,0}$.

\medskip

The result we will need reads~:
\begin{proposition}
Assume that $\delta \in (\frac{2-n}{2}, \frac{n-2}{2})$ and $\alpha \in (0,1)$ are 
fixed. There exists some constant $c >0$ and, for all $S >s_0$, there exists an operator 
\[
{\cal G}_S : {\cal C}^{0, \alpha}_\delta ([-S, S] \times S^{n-1}, {\mathfrak H}(n,m) 
) \longrightarrow  {\cal C}^{2, \alpha}_\delta ([-S, S] \times S^{n-1}, {\mathfrak H}
(n,m) ),
\]
such that, for all $f \in {\cal C}^{0, \alpha}_{\delta}([-S,S] \times S^{n-1}, 
{\mathfrak H}(n,m) )$, the function $w ={\cal G}_S (f)$ is a solution of 
\[
{\cal L} w  =  f 
\]
in $(-S,S) \times S^{n-1}$ with  $ w =0$ on $ \{\pm S\} \times S^{n-1}$. 
Furthermore, $\|w\|_{{\cal C}^{2, \alpha}_{\delta}} \leq c \, \|f\|_{{\cal C}^{0, 
\alpha}_{\delta}} $. 
\label{pr:3}
\end{proposition}
{\bf Proof :} Our problem being linear, we can assume without loss of generality that
\[
\sup_{[-S, S]\times S^{n-1}} \varphi^{- \delta} \, |f|  =1.
\]

Observe that, it follows from Proposition~\ref{pr:2} that, when restricted to the 
space of functions $w$ such that $w(s, \cdot)$ is orthogonal to $e_0, \ldots, 
e_n$ in the $L^2$ sense on $S^{n-1}$, the operator ${\cal L}$ is injective over 
$(-S, S)\times S^{n-1}$. Also, if $s >s_0$ then  ${\cal L}$ is injective over 
$(-S, S)\times S^{n-1}$ when restricted to functions which are even and only 
depend on $s$. As a consequence, for all $ S > s_0$, we are able to solve 
${\cal L} v  = f$, in $(-S,S)\times S^{n-1}$, with $v =0$ on $\{\pm S\}\times 
S^{n-1}$.

\medskip

We claim that there exists some constant $c>0$, independent of $S>s_0$ and of 
$f$, such that
\[
\sup_{[-S,S] \times S^{n-1}} | \varphi^{-\delta} \, w  |\leq c.
\]
Observe that the result is true when $S>1$ stays bounded. We argue by 
contradiction and assume that the result is not true. In this case, there would 
exist a sequence $S_k > 1$ tending to $+\infty$, a sequence of functions 
$f_{k}$ satisfying
\[
\sup_{[-S_k, S_k] \times S^{n-1}}  \varphi^{-\delta}  \, |f_k | = 1,
\]
and a sequence $v_k$ of solutions of ${\cal L} v_k = f_k$, in $(-S_k,S_k)
\times S^{n-1}$, with $v_k=0$ on $\{\pm S_k\}\times S^{n-1}$ such that 
\[
A_k :=  \sup_{[-S_k , S_k] \times S^{n-1}} | \varphi^{-\delta} \, v_k | 
\longrightarrow +\infty.
\]
Let us denote by $(s_k, \theta_k)\in [0, S_k)\times S^{n-1}$, a point where the 
above supremum is achieved, observe that all the functions we consider are 
even in the $s$ variable.  We claim that the sequence $S_k-s_k$ remains 
bounded away from $0$. Indeed,  since $v_k$ and $(\del_s^2 + 
\Delta_{S^{n-1}} ) \, v_k$ are both bounded by a constant (independent of 
$k$) times $\varphi^\delta (S_k) \, A_k$ in $[S_k -1 , S_k]\times S^{n-1}$ 
and since  $v_k=0$ on $\{S_k \}\times S^{n-1}$, we may apply standard 
elliptic estimates and conclude that the gradient of $v_k$ is also uniformly 
bounded by a constant times $\varphi^\delta (S_k) \, A_k$ in $[S_k- 
\frac{1}{2}, S_k ]\times S^{n-1}$. As a consequence the above supremum 
cannot be achieved at a point which is too close to $S_k$.  Therefore, up 
to some subsequence, we may also assume that the sequence $S_k -s_k$ 
converges to $S^*\in (0, +\infty]$.  We now distinguish a few cases 
according to be the behavior of the sequence $s_k$, which, up to a 
subsequence, can be assumed to converge in $[0, +\infty]$.

\medskip

\noindent
We define the sequence of rescaled functions
\[
\tilde{v}_k (s, \theta) : =  \frac{\varphi^{-\delta} (s_k)}{A_k} \, v_k (s+ 
s_k, \theta).
\]

\medskip

\noindent
{\em Case 1 :} Assume that the sequence $s_k$ converges to $s_*\in {\R}$. 
After the extraction of some subsequences, if this is necessary, we may assume
 that the sequence $\tilde{v}_k (\cdot -s_*, \cdot)$ converges on compact to 
$v$ some nontrivial solution of 
\[
{\cal L} v  =0 ,
\]
in ${\R} \times S^{n-1}$. Furthermore
\begin{equation}
\sup_{{\R} \times S^{n-1}} \varphi^{-\delta}  \, | v  | =\varphi^{-\delta} (s_*).
\label{eq:3-2}
\end{equation}
Moreover, for each $s \in {\R}$, the function $v (s, \cdot)$ is orthogonal in the 
$L^2$ sense to $e_1, \ldots, e_n$ on $S^{n-1}$. But, the result of 
Proposition~\ref{pr:2} together with the fact that $\Psi^{\pm,0}$ are the only 
solutions of ${\cal L} w=0$ which only depend on $s$, implies that $v \equiv 0$, 
contradicting (\ref{eq:3-2}).

 \medskip

\noindent
{\em Case 2} : Assume that the sequence $s_k$ converges to $+\infty$. After 
the extraction of some subsequences, if this is necessary, we may assume that 
the sequence $\tilde{v}_k$ converges to $v$ some nontrivial solution of 
\[
\Delta_0 \, v = 0,
\]
in $(- \infty, S^*)\times S^{n-1}$, with boundary condition $v =0$, if $S^*$ is finite. 
Furthermore
\begin{equation}
\sup_{(-\infty, S^*) \times S^{n-1}} e^{- \delta s} \, | v  | = 1.
\label{eq:3-3}
\end{equation}
independently of the fact that $S^*$ finite or is not, this case is easy to rule out 
using the eigenfunction decomposition of $v$ 
\[
v = \sum_{j \in {\mathfrak J}-\{0\}} v_j \, e_j.
\]
Indeed, $v_j$ has to be a linear combination of the functions $e^{\pm \gamma_j s}$
 (where $\gamma_j$ has been defined in (\ref{eq:3-4})) and is bounded by 
$e^{\delta s}$. Since we have assumed that $\delta \in (-\frac{n+2}{2},
 \frac{n+2}{2})$ and since $j \geq n+1$, it is easy to see that all $v_j \equiv 0$, 
contradicting (\ref{eq:3-3}). 

\medskip

We have reached a contradiction in all cases, hence, the proof of the claim is 
finished. To complete the proof of the Proposition, it suffices  sum the two results 
we have just obtained and apply Schauder's estimates in order to get the relevant 
estimates for all the derivatives. \hfill $\Box$

\medskip

We will also need some properties of the Poisson operator for $\Delta_0$ on 
$[0, \infty)\times S^{n-1}$. The result we will need is standard and a proof can 
be found, for example, in \cite{Fak-Pac}~: 
\begin{lemma}
There exists $c>0$ such that, for all $g \in  {\cal C}^{2,\alpha}(S^{n-1}, 
{\mathfrak H}(n,m)) $, there exists a unique $w \in 
{\mathcal C}^{2,\alpha}_{\frac{2-n}{2}}([0, +\infty) \times S^{n-1})$ solution 
of 
\begin{equation}
\left\{ 
\begin{array}{rllll} 
\Delta_0 w & =  & 0 \qquad & \mbox{\rm in}\quad (0,+\infty) \times 
S^{n-1}\\[3mm] 
         w & =  &  g \qquad & \mbox{\rm on} \quad \{ 0 \}\times S^{n-1}.
\end{array}  
\right. 
\label{eq:3-5}
\end{equation}
Furthermore,  we have $|| w ||_{{\cal C}^{2, \alpha}_{\frac{2-n}{2}} } 
\leq c \,  || g||_{{\cal C}^{2,\alpha}}$ and, for all $s >0$ the function $w(s, \cdot)$ 
is invariant with respect to the action of ${\mathfrak H}(n,m)$.
\label{le:1}
\end{lemma}
The idea behind the proof of this result is that one can use the eigenfunction 
decomposition of $g$ to obtain an explicite solution of (\ref{eq:3-5}) together 
with the estimate. In the remaining of the paper, we will denote by ${\cal P}(g)$ 
the solution of (\ref{eq:3-5}).

\subsection{The nonlinear problem}

We fix $\rho \in [0,1]$ and, for all $\e \in (0, \rho )$, we define $s_\e >0$ by the identity
\[
 \rho =  \e  \, \varphi (s_\e)   > 0 .
\]
Let us notice that, as $\e$ tends to $0$, we have  
\[
s_\e \sim  -   \log \e 
\]
In order to parameterize the unit $n$-catenoid we use  (\ref{eq:3-6}) and define the 
outer unit normal $N_0 $ as in  (\ref{eq:3-7}). Let us define a smooth function 
$\xi_\e : {\R} \longrightarrow [-1,1]$ which  satisfies $\xi_\e  = - 1 $ for $s \geq 
s_\e - 1$, $\xi_\e  = 1 $ for $s \leq 1 - s_\e$ and $\xi_\e  = - \ds 
\frac{\del_s \varphi}{\varphi}$ for $ |s| \leq s_\e -2$ and which interpolates 
smoothly between those two functions when $|s|\in [s_\e-2, s_\e-1]$. We 
consider the vector field 
\[
N_\e (s, \theta) : =  (\sqrt{1-\xi_\e^2 (s)} \, \theta , \xi_\e (s)).
\]
It turns out that this vector field is a perturbation of the unit normal $N_0$, 
and in fact, we have for all $k\geq 0$
\[
\left| \nabla^k \left( N_\e  \cdot N_0  - 1  \right) \right|\leq c_k \, \e^{2n-2},
\]
for all $|s| \geq s_\e - 2$. 

 \medskip

We now look for all minimal hypersurfaces close to the unit $n$-catenoid which 
has been rescaled by a factor $\e$. The hypersurfaces we are looking for will
be parameterized by 
\[
X_w : =  \e \, X_0 + w \, \varphi^{\frac{2-n}{2}} \,  N_\e ,
\]
for $(s, \theta) \in [- s_\e, s_\e] \times S^{n-1}$ and for some small function 
$w$.  It follows from (\ref{eq:3-1}) that such an hypersurface is minimal if 
and only if $w$ satisfies a nonlinear equation of the form
\begin{equation}
{\cal L} w = \bar{Q}_\e(w),
\label{eq:eqsd}
\end{equation}
where 
\[
\bar{Q}_\e(w)  := \ds  \e^{2n-2}\, L_\e w + \e \, \varphi^{\frac{2-n}{2}} \, 
\bar{Q}_{2,\e} \, \left(\varphi^{-\frac{n}{2}} \, \e^{-1} \, w \right) + \ds \e 
\, \varphi^{\frac{n}{2}} \,\bar{Q}_{3,\e} \,  \left(\varphi^{-\frac{n}{2}} \, 
\e^{-1} \, w \right). 
\]
Here $\bar{Q}_{2,\e}$ and $\bar{Q}_{3,\e}$ enjoy  properties which are 
similar to  those enjoyed by $Q_2$ and $Q_3$, namely (\ref{eq:useful1}) 
and (\ref{eq:useful2}) still hold uniformly in $\e \in (0,\rho)$.  The linear 
operator $ \e^{2n-2}\, L_\e$  represents the difference between the linearized 
mean curvature operator for hypersurfaces parameterized using the vector 
field $N_0$ and those parameterized using the vector field $N_\e$. The 
operator $L_\e$ has coefficients which are supported in $([-s_\e, 2 - s_\e] 
\cup [s_\e-2, s_\e])\times S^{n-1}$ and which are uniformly  bounded  
in ${\cal C}^{0, \alpha}$ topology. The details of the derivation of this 
formula can be found, for example, in \cite{Maz-Pac-Pol} or in 
\cite{Fak-Pac}.

\medskip

{\bf Solutions of (\ref{eq:eqsd}) which are parameterized by their 
boundary data~:} We  fix $\delta \in (\frac{2-n}{2}, \frac{n-2}{2})$, 
$\alpha \in (0, 1)$ and $\kappa >0$. Given $h \in {\cal C}^{2,\alpha}
(S^{n-1}, {\mathfrak H} (n,m))$ whose norm satisfies
\[
|| h  ||_{{\cal C}^{2,\alpha} } \leq  \kappa \, \e^{n-1} ,
\]
we set
\[
g : = \varphi^{\frac{n-2}{2}} (s_\e) \, h ,
\]
and we define
\begin{equation}
\tilde{w}  : = {\cal P}_{s_\e}(g) (s_\e - \cdot , \cdot) + 
{\cal P}_{s_\e}(g) (\cdot + s_\e , \cdot ) \in {\cal C}^{2, \alpha}
([-s_\e, s_\e]\times S^{n-1}, {\mathfrak H}(n,m)).
\label{eq:he}
\end{equation}
We know from Lemma~\ref{le:1} that
\begin{equation}
\| \tilde{w}\|_{{\cal C}^{2, \alpha}_{\frac{n-2}{2}}} \leq c\, 
\e^{\frac{n-2}{2}} \, \| g \|_{{\cal C}^{2, \alpha}}.
\label{eq:3-9}
\end{equation}

Now, if we write $w = \tilde{w} + v$, we wish to find a function $v \in 
{\mathcal C}^{2,\alpha}_{\delta}([-s_\e, s_\e]  \times S^{n-1}, 
{\mathfrak H}(n,m))$ such that 
\begin{equation}
\left\{ \begin{array}{rlll}
{\cal L} v & = &  \ds \bar{Q}_\e(\tilde{w}+v) - {\cal L} \tilde w 
\qquad & \mbox{in} \qquad (- s_\e, s_\e )\times S^{n-1} \\[3mm]
             v & = & 0 \qquad & \mbox{on} \qquad \{ \pm s_\e \} \times 
S^{n-1}.
\end{array}
\right.
\label{eq:10p}
\end{equation}
To obtain a solution of this equation, it is enough to find a fixed point 
of the mapping
\[
{\cal N}_\e (v) :=  {\cal G}_{s_\e }\left(\bar{Q}_\e (\tilde{w} + v) - 
{\cal L}\tilde w\right),
\]
where the operator ${\cal G}_{s_\e}$ has been defined in 
Proposition~\ref{pr:1}. Using (\ref{eq:3-9}) together with 
Proposition~\ref{pr:1} and the properties of $\bar Q_\e$, we can 
estimate 
\[
\| \e^{2n-2} \, L_\e \tilde w  - {\cal L} \tilde w \|_{{\cal C}^{0, 
\alpha}_{\delta}} \leq c \, \left( 1  + \e^{\frac{3n-2}{2} + \delta}  
\right) \, \|h \|_{{\cal C}^{2, \alpha}},
\]
\[
\| \e  \, \varphi^{\frac{2-n}{2}} \, \bar{Q}_{2,\e} \, 
\left(\varphi^{-\frac{n}{2}} \, \e^{-1} \, \tilde w \right) 
\|_{{\cal C}^{0, \alpha}_{\delta}} \leq c \, \left( \e^{-1} + 
\e^{\frac{n}{2}+\delta} \right) \, \|h \|_{{\cal C}^{2, \alpha}}^2
\]
and finally, there exists $\e_0>0$ (which depends on $\kappa$) such that
for all $\e \in (0, \e_0)$ we have
\[
\| \e  \, \varphi^{\frac{n}{2}} \,\bar{Q}_{3,\e} \,  
\left(\varphi^{-\frac{n}{2}} \, \e^{-1} \,  \tilde w \right)\|_{
{\cal C}^{0, \alpha}_{\delta}}\leq c \,  (\e^{-2} + 
\e^{\frac{2-n}{2}+\delta} )  \, \| h \|_{{\cal C}^{2, \alpha}}^3 ,
\]
In the above estimates, the constant $c>0$ does not depend on $\e$, 
nor on $\kappa$. Observe that in order to obtain the last estimate, 
we have implicitely used that fact that $\|h\|_{{\cal C}^{2, \alpha}}$ 
is small enough so that we can apply (\ref{eq:useful2}), or rather its 
counterpart for $\bar Q_{3, \e}$.  This explains 
the restriction of $\e \in (0, \e_0)$ which is needed.

\medskip

It is then a simple exercise to show that for any fixed $\kappa >0$, 
there exist $c >0$ and $\e_0 >0$, such that, for all $\e \in (0, 
\e_0)$, the nonlinear mapping ${\cal N}_\e$ is a contraction in the 
ball of radius 
\[
R (\e, h) : =  c  \, \|h\|_{{\cal C}^{2, \alpha}},
\]
in ${\cal C}^{2, \alpha}_\delta ([-s_\e, s_\e]  \times S^{n-1}, 
{\mathfrak H}(n,m))$ into itself, and hence ${\cal N}_\e$ has a 
unique fixed point $v_h$ in this ball. Therefore, the function $\tilde 
w + v_h$ is a solution of (\ref{eq:eqsd}) whose boundary data 
is given by $h$, up to a constant function. We can even choose the 
constant  $c$ to be independent of $\kappa$, but this will not be
useful.

\medskip

{\bf Family of minimal hypersurfaces close to $n$-catenoid~:} We 
summarize the results we have obtained so far and translate them in 
the geometric framework. Let us fix $\delta \in (\frac{2-n}{2}, 
\frac{n-2}{2})$, $\alpha \in (0, 1)$ and $\kappa >0$. There
 exists $c >0$ and $\e_0 >0$ such that for all
 $\e \in (0, \e_0)$ and for all $h \in  {\cal C}^{2, \alpha} (S^{n-1}, 
{\mathfrak H}(n,m)))$ satisfying 
\[
|| h  ||_{2,\alpha} \leq  \kappa \, \e^{n-1},
\]
there exists a  minimal hypersurface, which will be denoted by 
$C_\e (h) \subset {\R}^{n+1}$, and which is parameterized by 
\[
X_w \: = \e \, X_0 + w \, \varphi^{\frac{2-n}{2}} \, N_\e \qquad
 \mbox{in}\qquad  [- s_\e, s_\e] \times S^{n-1},
\]
for some function $w$ satisfying 
\[
\|w\|_{{\cal C}^{2, \alpha}_{\frac{2-n}{2}}} \leq c \, 
\|h\|_{{\cal C}^{2, \alpha}}.
\]

\medskip

This hypersurface is symmetric with respect to the hyperplane 
$z=0$ and further inherits all the symmetries  induces by the 
symmetries used to define the function spaces in 
Definition~\ref{de:1}, hence it is invariant with respect to the 
action of $O(n-m, {\R}) \otimes {\mathfrak D}_m \otimes  
\{\pm I_1\}\subset O(n+1, {\R})$. Furthermore, if we perform 
the change of variable 
\[
r : = \e \, \varphi (s),
\]
we see that near its upper boundary, this hypersurface is the 
graph of the function
\[
x \in \overline{B^n_{\rho}} \setminus B^n_{\rho/2} 
\longrightarrow \e \, c_\infty - W_{h}(x) - V_{\e, h}(x),
\]
over the $z=0$ hyperplane. Here $W_h$ denotes the (unique) 
harmonic extension of the boundary data $h$ in $B^n_{\rho}$ 
and the function $V_{\e, h} $ satisfies
\[
\|V_{\e, h} \|_{{\cal C}^{2, \alpha}} \leq  c_0 \, \e^{n-1}.
\] 
for some constant $c_0$ which does not depend on $\kappa$ 
nor on $\e$. Here the norms are taken over $\overline{B^n_{\rho}}-
B^n_{\rho/2}$.  This last claim, which is a key point of our analysis,
follows from (\ref{eq:dldl}). Indeed, when $h=0$, $C_\e (0)$
 is just a rescaled $n$-catenoid and, using (\ref{eq:dldl}) we 
see that its upper end is the graph of the function
\[
x \longrightarrow \e \,  c_\infty  +  {\cal O} ( \e^{n-1}  \, r^{2-n}).
\]
We have also used the fact that the solution of (\ref{eq:10p}) we 
have constructed is equal to $\tilde w + v_h$ where $\tilde w$, 
defined in (\ref{eq:he}), is linear in $h$ and where $v_h$ can 
be estimated by a constant (independent of $\e$ and $\kappa$) 
times $\|h\|_{{\cal C}^{2, \alpha}} \,  \varphi^{\delta}$.
Essentially the constant $c_0$ arises from the term ${\cal O} 
( \e^{n-1}  \, r^{2-n})$ in the above expansion, the contributions
of $v_h$ and the pertubation caused by the change of variable 
being negligeable when $\e$ is chosen small enough.

\medskip

Observe that, reducing $\e_0$ if this is necessary, we can 
assume that the mapping $h \longrightarrow V_{\e, h}$ is 
continuous and in fact smooth. With little work we also find 
that
\begin{equation}
\|V_{\e, h_2} -V_{\e, h_1}\|_{{\cal C}^{2, \alpha}}\leq c\, 
\e^{\frac{n-2}{2}-\delta} \, \|h_2-h_1\|_{{\cal C}^{2, \alpha}} ,
\label{eq:7-1}
\end{equation}
for some constant $c >0$ which does not depend on $\e$. The 
norm on the left hand side of this inequality is 
understood to be the norm on $\overline{B^n_{\rho}}-
B^n_{\rho/2}$. Again the constant $c$ can be chosen to be 
independent of $\kappa$ but this will be irrelevent for the 
remaining of the analysis. 

\section{Minimal hypersurfaces which are graphs over an 
hyperplane}

We are now concerned with both the mean curvature and 
the linearized mean curvature operator for hypersurfaces 
which are graphs over the $z=0$ hyperplane, in ${\R}^{n-m} 
\times T^m \times {\R}$.

\subsection{The mean curvature operator for graphs}

We assume that $n\geq 3$ and $1 \leq m \leq n-1$ are fixed. 
Further assume that $T^m \in {\cal T}^m$ is fixed. Then, for 
any function  $u $, defined in  ${\R}^{n-m} \times T^m $, 
which is at least of class ${\cal C}^2$, we can define an 
hypersurface which is the graph of $u$
\[
{\R}^{n-m} \times T^m  \ni (x_1,x_2)  \longrightarrow (x_1,
x_2, u(x_1,x_2)) \in {\R}^{n-m} \times T^m \times {\R}.
\]
Recall that the mean curvature of this hypersurface, with 
downward pointing unit normal, is then given by
\begin{equation}
H_u : =  - \frac{1}{n} \mbox{div} \, \left( \frac{\nabla u}{(1
+|\nabla u|^2)^{1/2}}\right) .
\label{eq:nbnb}
\end{equation}

\subsection{Linear analysis}

We define the function spaces which are adapted to the analysis 
of the Laplacian in $T^m \times {\R}^{n-m}$. Our main concern 
will be the asymptotic behavior of the functions as $|x_1|$ tends 
to $+\infty$.
\begin{definition} 
For all $k \in {\N}$,  $\alpha \in (0,1)$ and $\nu \in {\R}$, the 
space ${\cal C}^{k,\alpha}_\nu ( {\R}^{n-m} \times T^m) $ is 
defined to be the space of functions $ w \in {\cal C}^{k,
\alpha}_{\mathrm loc} ( {\R}^{n-m} \times T^m)$ for which 
the following norm is finite
\[
\|w \|_{{\cal C}^{k, \alpha}_\nu} : =  |w |_{{\cal C}^{k, \alpha } 
( B^{n-m}_{1} \times T^m) } + \sup_{r > 1/2} r^{-\nu} \, | w (r \, 
\cdot) |_{{\cal C}^{k,\alpha}(  (B^{n-m}_2-B^{n-m}_1) \times 
r^{-1} T^m)} .
\] 
Here $| \, \, \, |_{{\cal C}^{k, \alpha}( \Omega)}$ denotes the 
H\"older norm in $\Omega$.
\label{de:22}
\end{definition}
To get a better undertanding of these weighted spaces, if $T_m 
={\R}^m / A {\mathbb Z}^m$, we can identify any function defined 
on  $ {\R}^{n-m} \times T^m$ with a function defined on 
${\R}^{n-m}\times {\R}^m$ which has $\{0\} \otimes A \, {\mathbb Z}^m$ as 
its group of periods. In which case functions which belong to 
${\cal C}^{k,\alpha}_\nu ( {\R}^{n-m} \times T^m) $ are identified 
with functions defined on ${\R}^{n-m}\times R^m$, which are 
bounded by a constant times $(1+|x_1|)^\nu$, whose first derivative
 is bounded by a constant times $(1+|x_1|)^{\nu-1}$ (if $k\geq 1$), 
and so on.

\medskip

As in the previous section, we will only work with functions having 
some special symmetry. Therefore, we introduce the~:
\begin{definition} 
For all $k \in {\N}$,  $\alpha \in (0,1)$ and $\nu \in {\R}$, the space 
${\cal C}^{k,\alpha}_\nu ( {\R}^{n-m} \times T^m, {\mathfrak H}
(n,m)) $ is defined to be the space of functions $ w \in 
{\cal C}^{k,\alpha}_{\mathrm loc} ( {\R}^{n-m} \times T^m)$ 
which are invariant under the action of ${\mathfrak H}(n,m)$.
\label{de:2}
\end{definition}
Observe that, because of the invariance of our function space 
with respect to the action of ${\mathfrak H} (n,m)$, any function 
$w \in {\cal C}^{k,\alpha}_\nu ( {\R}^{n-m} \times T^m,  
{\mathfrak H}(n,m))$ can be decomposed as
\[
w (x_1,x_2)  = \sum_{i \in {\mathfrak I}} w_i (r_1)  \, E_i (x_2),
\]
where ${\mathfrak I}(m) \subset {\mathbb N}$ has been defined 
in (\ref{eq:2-1}) and where 
\[
r_1 : = |x_1|.
\]

\medskip

To begin with let us treat the easy case where $ 1 \leq m \leq n-3$. We 
have the~: 
\begin{proposition}
Assume that $ 1 \leq m \leq n-3$. Given $\nu \in (2+m-n, 0)$ and 
$\alpha \in (0,1)$. There exist some constant $c >0$ and  an operator 
\[
G : {\cal C}^{0, \alpha}_{\nu -2} ({\R}^{n-m} \times T^m , 
{\mathfrak H}(n,m)) \longrightarrow  {\cal C}^{2, \alpha}_\nu 
({\R}^{n-m} \times T^m,  {\mathfrak H}(n,m)) ,
\]
such that, for all $f \in {\cal C}^{0, \alpha}_{\nu -2}({\R}^{n-m} 
\times T^m,  {\mathfrak H}(n,m))$, the function $w =G (f)$ is a 
solution of 
\[
\Delta w = f ,
\]
in ${\R}^{n-m} \times T^m$. Furthermore, 
$\|w\|_{{\cal C}^{2, \alpha}_\nu} \leq c \, 
\|f\|_{{\cal C}^{0, \alpha}_{\nu-2}} $. 
\label{pr:41}
\end{proposition}
{\bf Proof~:} The proof of the result is simplified by the fact that 
\[
\Delta |x_1|^\nu  =  - \nu \, ( n -m-2 -\nu) \, |x_1|^{\nu-2}
\]
Hence, the function $w (x_1, x_2) : = |x_1|^{\nu}$, which is defined 
in $({\R}^{n-m} -\{0\} )\times T^m$ can be used as a barrier function 
to prove, for any 
$f \in {\cal C}^{0, \alpha}_{\nu -2}({\R}^{n-m} \times T^m,  
{\mathfrak H}(n,m))$, the existence of a solution of
\[
\Delta w = f ,
\]
in ${\R}^{n-m} \times T^m$. Furthermore, it also yields the estimate
\[
|w (x_1, x_2) | \leq c \, \|f\|_{{\cal C}^{2 , \alpha}_{\nu-2}}\, |x_1|^\nu
\]
for some constant which does not depend on $f$. The maximum 
principle then implies that
\[
|w (x_1, x_2) | \leq c \, \|f\|_{{\cal C}^{2 , \alpha}_{\nu-2}}\, (1+ |x_1|)^\nu .
\]
 Starting from this, Schauder's estimates yield
\[
\|w  \|_{{\cal C}^{2, \alpha}_\nu} \leq c \, \|f\|_{{\cal C}^{2 , \alpha}_{\nu-2}}.
\]
Details are left to the reader. \hfill $\Box$

\medskip

When $m=n-2$ or $m=n-1$, the previous result has to be modified 
since $2+m-n \geq 0$ in these two cases . To this aim, 
we choose  $\chi$ a cutoff function defined on ${\R}$ such that 
$\chi \equiv 1$ for $t\geq 2$ and $\chi \equiv 0$ when $t \leq 1$. When $m=n-2$, 
we define the space 
\[
{\cal D}_{2} : = \mbox{Span}\{ \chi (r_1) \, \log  r_1\} \subset 
{\cal C}^{\infty}({\R}^2),
\]
and when $m=n-1$, we set 
\[
{\cal D}_{1} : = \mbox{Span}\{\chi(r_1) \, r_1 \} \subset {\cal C}^{\infty}(\R) .
\]
This time we have the~: 
\begin{proposition}
Assume that $m=n-2$ or $m=n-1$. Given $\nu \in (-\infty, 0)$ and $\alpha \in (0,1)$. 
There exist some constant $c >0$ and  an operator 
\[
G : {\cal C}^{0, \alpha}_{\nu -2} ({\R}^{n-m} \times T^m , {\mathfrak H}(n,m)) 
\longrightarrow  {\cal C}^{2, \alpha}_\nu ({\R}^{n-m} \times T^m,  
{\mathfrak H}(n,m)) \oplus {\cal D}_{n-m},
\]
such that, for all $f \in {\cal C}^{0, \alpha}_{\nu -2}({\R}^{n-m} \times T^m,  
{\mathfrak H}(n,m))$, the function $w =G (f)$ is a solution of 
\[
\Delta w = f ,
\]
in ${\R}^{n-m} \times T^m$. Furthermore, $\|w\|_{{\cal C}^{2, \alpha}_\nu 
\oplus {\cal D}_{n-m}} \leq c \, \|f\|_{{\cal C}^{0, \alpha}_{\nu-2}} $. 
\label{pr:4}
\end{proposition}
{\bf Proof~:} We decompose
\[
f  = f_0 + \sum_{i \in {\mathfrak I}-\{0\}} f_i  \, E_i ,
\]
and adopt the notation $f = f_0 + f'$. We look for a solution $w$ which will 
also be decomposed as
\[
w =w_0 + \sum_{i \in {\mathfrak I}-\{0\}} w_i \, E_i,
\]
and again we set $w = w_0 + w'$. For notational convenience,  $f', v', w', 
\ldots$ will denote functions whose eigenfunction decomposition only involves 
indices $i \in {\mathfrak I} (m) -\{0\}$.

 \medskip

Observe that, because of the invariance of our problem with respect to the 
action of ${\mathfrak H}_n$,  the Laplacian in ${\R}^{n-m} \times T^m$ 
reduces to the study of the operator
\[
L : = \del_{r_1}^2 + \frac{n-m-1}{r_1} \, \del_{r_1} + \Delta_{T^m},
\]
where we have set $r_1 : = |x_1|$.

\medskip

{\bf Step 1 :} We would like to prove the existence of $w'$ and also obtain 
the relevant estimate.  Our problem being linear, we may always assume that 
\[
\sup_{{\R}^{n-m} \times T^m }  (1+r_1)^{2-\nu}  \,  | f' | = 1 .
\]
Obviously $\Delta$, or $L$,  is injective over any  $B^{n-m}_R \times T^m$. 
As a consequence, for any $R > 1$ we are able to solve $\Delta v' = f'$, in 
$B^{n-m}_R \times T^m$, with $v' =0$ on $\del B^{n-m}_R \times T^m$.

\medskip

We claim that, there exists a constant $c>0$, independent of $R>1$ and 
of $f'$, such that
\[
\sup_{B^{n-m}_R \times T^m}  (1+r_1)^{-\nu} \, | v' |\leq c.
\]
Observe that the result is certainly true if we assume that $R$ remains 
bounded. We argue by contradiction and assume that the claim is not true. 
In this case, there would exist a sequence $R_k > 1$ tending to $+\infty$, 
a sequence of functions $f_k'$ satisfying
\[
\sup_{B^{n-m}_{R_k} \times T^m} (1+r_1)^{2-\nu}  \, | f'_k | = 1,
\]
and a sequence $v'_k$ of solutions of ${\cal L} v'_k = f_k'$, in 
$B^{n-m}_{R_k} \times T^m$, with $v'_k=0$ on $\del B^{n-m}_{R_k} 
\times T^n$, such that 
\[
A_k :=  \sup_{B^{n-m}_{R_k} \times T^m} (1+r_1)^{-\nu} \, | v_k' | 
\longrightarrow +\infty.
\]
Let us denote by $(x_{1,k}, x_{2,k})\in  B^{n-m}_{R_k} \times T^m$, a 
point where the above supremum is achieved. We now distinguish a few 
cases according to the behavior of the sequence $r_{1,k} :=  |x_{1,k}|$ 
which, up to a subsequence can always be assumed to converge in 
$[0, +\infty]$. Observe  that, as in the proof of Proposition~\ref{pr:1},
 the sequence $R_k-r_{1,k}$ remains bounded away from $0$. 

\medskip

\noindent
We define the sequence of rescaled functions
\[
\tilde{v}'_k (x_1, x_2) : =  \frac{(1+r_{1,k})^{-\nu}}{A_k} 
\, v'_k (r_{1,k} \, x_1, r_{1,k} \, x_2).
\]

\medskip

\noindent
{\em Case 1 :} Assume that the sequence $r_{2,i}$ converges to
 $r_{2,\star} \in [0, \infty)$. After the extraction of some subsequences, 
if this is necessary, we may assume that the sequence $\tilde{v}_k' 
(\cdot / r_{1,k} , \cdot / r_{1,k})$ converges to some nontrivial solution of 
\begin{equation}
L v' =0 ,
\label{eq:4-1}
\end{equation}
in ${\R}^{n-m}\times T^m$. Furthermore
\begin{equation}
\sup_{{\R}^{n-m} \times T^m} (1+r_1)^{- \nu}  \, |v' | = 
(1+ r_{1,\star})^{-\nu}.
\label{eq:4-2}
\end{equation}
But the maximum principle implies that $v$ is identically equal to $0$. 
This clearly contradicts (\ref{eq:4-2}).

 \medskip

\noindent
{\em Case 2} : Assume that the sequence $r_{2,k}$ converges to $+\infty$. 
After the extraction of some subsequences, if this is necessary, we may 
assume that the sequence $\tilde{v}_k'$ converges to some nontrivial solution of 
\[
\del_{r_1}^2 v' + \frac{n-m -1}{r_1} \, \del_{r_1} v'  + 
\Delta_{{\R}^{n-m}} v' = 0,
\]
in $({\R}^{n-m}-\{0\}) \times {\R}^{m}$, which does not depend on $x_2$.  
This last claim follows from the fact that the functions $x_2 \rightarrow 
\tilde v_k'(x_1, x_2)$ have a group of period given by $r_{1,k}^{-1} 
\, A \, {\mathbb Z}^m$, if $T^m = {\R}^m / A \, {\mathbb Z}^m$. In 
addition, $|\nabla_{x_2} \tilde v_k'$
is bounded by a constant depending on $x_1$. PAssing to the limit, we
see that $v'$ does not depend on $x_2$. 

\medskip

Furthermore
\begin{equation}
\sup_{({\R}^{n-m}-\{0\}) \times R^m} r_1^{-\nu} \, | v' | = 1.
\label{eq:4-3}
\end{equation}
It should be clear that $v' \equiv 0$, contradicting (\ref{eq:4-3}). 

\medskip

Since we have obtained a contradiction in both cases, this 
finishes the proof of the claim.

\medskip

{\bf Step 2 :}  We now turn our attention to the existence of $w_0$ 
as well as the relevant estimate. Again, our problem reduces to one 
ordinary differential equation since we now have to solve 
\[
\del_{r_1}^2 w_0 + \frac{n-m-1}{r_1} \del_{r_1} w_0  = f_0 .
\]
It is easy so check that $w_0$ is given by the formula
\[
w_0 =   \int_0^{r_2} \zeta^{1 + m -n}\, \int_0^\xi  t^{n - m -1} \, f_0 \, dt \, d\zeta - \ds\int_0^\infty \zeta^{1+m-n}\, \int_0^\xi  t^{n-m -1} \, f_0 \, dt \, d\zeta  .
\]  
In order to simplify the exposition, we will restrict our attention to the case where $m \neq n-2$ since, obvious modifications have to be done to treat the case $m = n - 2$. Granted the above formula, one can directly check that we can decompose, for all $r_2 >1$, $w_0 : = a_0 \, r_1^{2 + m-n} +\tilde{w_0}$, where 
\[
\begin{array}{rllll}
a_0 & := & \ds \frac{1}{n -m -2} \,  \int_0^\infty   t^{n-m -1} \, f_0 \, dt \, d\zeta , \\[3mm]
\tilde w_0 & : = & \ds \int_{r_2}^\infty \zeta^{1 + m - n}\, \int_\xi^\infty  t^{n-m -1} \, f_0 \, dt \, d\zeta .
 \end{array}
\]
Moreover, we have 
\[
|a_0| + \sup_{(0, \infty)} (1+ r_1)^{-\nu} \, | \tilde w_0| \leq c \, \sup_{(0, \infty)} (1+ r_1)^{2-\nu} \, |f_0| .
\]

\medskip

To complete the proof of the Proposition, it suffices to sum the two 
results we have just obtained and apply Schauder's estimates in order 
to get the relevant estimates for all the derivatives. \hfill $\Box$

\medskip

If $\rho >0$ is fixed, we define the space ${\cal C}^{k, \alpha}_{\nu}
({\R}^{n-m}\times T^m -B^n_\rho, {\mathfrak H}(n,m))$ as the 
space of restriction of functions of  ${\cal C}^{2, \alpha}_{\mu} 
({\R}^{n-m}\times T^m , {\mathfrak H}(n,m))$ to ${\R}^{n-m}
\times T^m -B^n_\rho$. This space is naturally endowed with the 
induced norm. 

\medskip

In order to simplify notations, we set 
\[
{\cal E}^{2, \alpha}_\nu : = {\cal C}^{2, \alpha}_{\nu}
({\R}^{n-m}\times T^m -B^n_\rho, {\mathfrak H}(n,m))
\]
when $1 \leq m \leq n-3$ and 
\[
{\cal E}^{2, \alpha}_\nu : = {\cal C}^{2, \alpha}_{\nu}
({\R}^{n-m}\times T^m -B^n_\rho, {\mathfrak H}(n,m)) 
\oplus {\cal D}_{n-m}
\]
when $m=n-2$ or $m=n-1$. We also define
\[
{\cal F}^{0, \alpha}_{\nu-2} : = {\cal C}^{0, \alpha}_{\nu}
({\R}^{n-m}\times T^m -B^n_\rho, {\mathfrak H}(n,m)) 
\]
when $1\leq m \leq n-1$. Using the previous result together with 
a standard perturbation result, we obtain the~:
\begin{proposition}
Assume that  $\nu \in (2+m-n,0)$ when $1 \leq m \leq n-3$, 
$\nu \in (-\infty, 0)$ when $m=n-2$ or $m=n-1$, and $\alpha \in (0,1)$ are fixed. 
There exist $\rho_0 >0$, $c >0$ and,  for all  for all $\rho \in (0, \rho_0)$, 
there exists an operator 
\[
G_\rho : {\cal F}^{0, \alpha}_{\nu -2}  \longrightarrow  {\cal E}^{2, \alpha}_{\nu} 
\]
such that, for all $f \in {\cal F}^{0, \alpha}_{\nu -2 }$, the function $w =G_\rho (f)$ 
is a solution of 
\[
\Delta w = f,
\]
in ${\R}^{n-m}\times T^m -B^n_\rho$, with $w \in {\R}$ on $\del 
B^n_\rho$. Furthermore, $\|w\|_{{\cal E}^{2, \alpha}_{\nu }} \leq 
c \, \|f\|_{{\cal F}^{0, \alpha}_{\nu-2}} $.  
\label{pr:5}
\end{proposition}

\subsection{The nonlinear problem}

Using (\ref{eq:nbnb}), one can check that the hypersurface parameterized by 
\[
{\R}^{n-m} \times T^m -B_\rho \ni (x_1,x_2) \longrightarrow (x_1, x_2, 
u(x_1,x_2)) \in {\R}^{n-m} \times T^m \times {\R},
\]
has mean curvature $0$ if and only if the function $u$ is a solution of 
\begin{equation}
\Delta u  - Q(u) = 0 ,
\label{eq:aaa}
\end{equation}
where we have set
\[
Q (u)  : = - \frac{1}{1+|\nabla u|^2} \, \nabla^2u \, ( \nabla u, \nabla u). 
\]

{\bf Solutions of (\ref{eq:aaa}) which are parameterized by their boundary 
data~:}  Let us assume that 
\[
\nu \in (2+m-n, 0) \qquad \mbox{when} \qquad 1\leq m \leq n-3,
\]
or that 
\[
\nu\in (-2, 0) \qquad \mbox{when} \qquad m=n - 2 
\]
or that 
\[
\nu \in (-\infty, 0)\qquad \mbox{when} \qquad m=n-1
\]
is fixed. The new restriction on $\nu$ when $m=n-2$ is needed to ensure
 that the nonlinear operator 
\[
u \longrightarrow \Delta u - Q(u) ,
\]
maps ${\cal E}^{2, \alpha}_{\nu}$ into ${\cal F}^{0, \alpha}_{\nu -2}$. 
Thanks to the result of Proposition~\ref{pr:5}, it is possible to apply the 
implicit function theorem to solve (\ref{eq:aaa}) with $w$ on 
$\del B^n_\rho$ equal to some given function $h  \in  {\cal C}^{2, \alpha}
(S^{n-1}, {\mathfrak H} (n,m))$ satisfying $\|h\|_{2, \alpha}\leq c_0$ 
for some fixed constant $c_0 >0$.  The solution of (\ref{eq:aaa}) provided
 by the implicit function theorem will  be denoted  by $w_h$. By construction, 
the graph of $w_{h}$ is a minimal hypersurface whose boundary is 
parameterized by the boundary data $h$.

\medskip

{\bf Family of minimal hypersurfaces which are close to ${\R}^{n-m} \times 
T^m$~} Let us summarize what we have proved.  We fix $\nu$ according 
to the above choice, $\alpha \in (0, 1)$ and $\kappa >0$. There exists 
$\e_0 >0$ and for all $\e \in (0, \e_0)$,  for all $h \in {\cal C}^{2, \alpha} 
(S^{n-1}, {\mathfrak H} (n,m))$ satisfying  
\[
\|h\|_{2, \alpha}\leq \kappa \, \e^{n-1}, 
\]
we have been  able to find a minimal hypersurface, which is a graph over 
$\Omega_\e$. This hypersurface, once translated by $\e\, c_\infty$ along 
the $z$-axis, will be denoted by $M_\e (h)$.

\medskip

Moreover, there exists a constant $c_{h,t}$ such that $M_\e (h)$ 
asymptotic to 
\[
\{(x_1, x_2,z) \in {\R}^{n-m} \times T^{m} \times  {\R} \quad | 
\quad z =  c_{h} \, \zeta_m (x_1)  \},
\]
where $\zeta_{n-1}(y) := |y|$, $\zeta_{n-2}(y) := \log |y|$ and 
$\zeta_m (y) := 0$, when $1\leq m \leq n-3$.  

\medskip

Also observe that the hypersurface $M_\e (h)$ inherits all the symmetries  
induces by the symmetries used to define the function spaces in 
Definition~\ref{de:2}, hence it is invariant with respect to the action of 
$O(n-m, {\R}) \otimes {\mathfrak D}_m \otimes \{I_1\}\subset 
O(n+1, {\R})$. 

\medskip

Furthermore, near its boundary this hypersurface can be parameterized 
as the graph of
\[
\overline{B^n_{2\rho}} - B^n_{\rho}\ni x  \longrightarrow  \e \, c_\infty 
 - \widehat W_h (x) - \widehat V_{h}(x),
\]
where $\widehat W_h$ is the unique (bounded) harmonic extension of the 
boundary data $h$ in ${\R}^{n-m} \times T^m - B^n_{\rho}$ which 
belongs to ${\cal E}^{2, \alpha}_{\nu}$. Here the function $\widehat V_{h}$ 
satisfies 
\[
\| \widehat V_{h}\|_{{\cal C}^{2, \alpha}}\leq c_0 \, \e^{2n-2},
\]
for some constant $c_0$ which depends on $\kappa$ but does not depend 
on $\e$. Here the norm is taken over $B^n_{2\rho}-B^n_{\rho}$.

\medskip

Observe that, reducing $\e_0$ if this is necessary, we can assume that the 
mapping $h \rightarrow \widehat V_{h}$ is continuous and in fact smooth. 
It follows from standard properties of the solutions obtained through the 
application of the implicit function theorem that
\begin{equation}
\| \widehat V_{h_2} - \widehat V_{h_1}\|_{{\cal C}^{2, \alpha}}\leq 
c\,  \e^{n-1} \, \|h_2-h_1\|_{{\cal C}^{2, \alpha}} ,
\label{eq:7-2}
\end{equation}
for some constant $c >0$ which does not depend on $\e$, but depends on $\kappa$. 
Here the norms are understood on $B^n_{2\rho}-B^n_\rho$.

\section{The gluing procedure}

We fix $\kappa >0$ large enough and apply the results of the previous sections. 
There exists $\e_0 >0$ and for all $g,h,  \in {\cal C}^{2, \alpha}(S^{n-1}, 
{\mathfrak H} (n,m))$  satisfying $\|g\|_{2, \alpha} \leq \kappa \, \e^{n-1}$ 
and $\|h\|_{2, \alpha} \leq \kappa \, \e^{n-1}$, we define the hypersurface 
$M_\e (g )$ and the hypersurface $C_\e (h )$. Our aim will be now to find 
$g$ and $h$  in such a way that 
\[
\left( M_\e ( g ) \cup   C_\e  (h) \right) \, \cap \, {\R}^{n-m}\times T^m 
\times (0, +\infty) ,
\]
is a ${\cal C}^{1}$ hypersurface. Then applying a reflection with respect to 
the hyperplane $z=0$, we will obtain a complete ${\cal C}^1$ hypersurface 
of ${\R}^{n-m} \times T^m \times {\R}$. Finally, it will remain to apply 
standard regularity theory to show that this hypersurface is in fact 
${\cal C}^\infty$.
 
\medskip

By construction, the two hypersurfaces  $M_\e (g, h_0)$ and $C_\e (h')$ are 
graphs over the $z=0$ hyperplane near their common boundary of the function 
\[
x \in B^n_{2 \rho}\setminus B^n_{\rho} \longrightarrow \e \, c_\infty 
- \widehat W_{g}(x) - \widehat V_{g}(x),
\]
for $M_\e ( g )$ and of the function
\[
x \in B^n_{\rho}\setminus B^n_{\rho/2} \longrightarrow \e \, c_\infty  - 
W_{h}(x) - V_{\e, h}(x),
\]
for $C_\e ( h )$.

\medskip

Hence, to produce a ${\cal C}^1$ hypersurface, it remains to solve the equations 
\begin{equation}
\left\{
\begin{array}{rllll}
 \widehat W_{g} + \widehat V_{g} & = & \ds  W_{h} - V_{\e, h} \\[3mm]
\ds  \del_r \widehat W_{g}  + \del_r \widehat V_{\e, g }& = & \ds 
 \del_r W_{h}+  \del_r V_{\e, h} ,
\end{array}
\right.
\label{eq:nn}
\end{equation}
where all functions are evaluated on $\del B^n_{\rho}$. The first identity is 
obtained by asking that the Dirichlet data of the two graphs on $\del B^n_{\rho}$ 
coincide and already ensures that the hypersurface is ${\cal C}^0$, while the second 
is obtained by asking that the Neumann data of the two graphs on $\del B^n_{\rho}$ 
coincide and ensures that the hypersurface will be of class ${\cal C}^1$.

\medskip

To this aim,, let us recall that the mapping ${\cal U}$
\[
{\cal U} :  h \in {\cal C}^{2, \alpha}(S^{n-1}, {\mathfrak H}(n,m)) \longrightarrow 
 \rho \, \del_r (W_h - \widehat W_h) (\rho \, \cdot ) \in {\cal C}^{1, \alpha} (S^{n-1} , 
{\mathfrak H}(n,m)),
\]
is an isomorphism. Indeed, this mapping is a linear first order elliptic pseudo-differential 
operator and, in order to check that it is an 
isomorphism, it is enough to prove that it is injective. Now if we assume that 
${\cal U} (h)=0$ then the function $w$ defined by $w := \widehat W_h$ in 
${\R}^{n-m}\times T^m -B^n_\rho$ and $w: =  W_h$ in $B^n_\rho$ is a global 
solution of $\Delta w =0$ in ${\R}^{n-m}\times T^m$, and furthermore, $w$ 
belongs to ${\cal E}^{2, \alpha}_{\nu} $. It is easy to check that necessarily $w 
\equiv 0$ and, as a consequence, $h\equiv 0$.  

\medskip

Using the above claim,  it is easy to see that (\ref{eq:nn}) reduces to a fixed 
point problem 
\[
(g, h ) = {\bf C}_\e (g, h ),
\]
in ${\cal E} : =  ({\cal C}^{2, \alpha}(S^{n-1}) , {\mathfrak H}(n,m))^2$. 
However,  (\ref{eq:7-1}) and (\ref{eq:7-2}) imply that ${\bf C}_\e :  {\cal E} 
\longrightarrow {\cal E}$ is a contraction mapping defined in the ball of radius 
$\kappa \, \e^{n-1}$ of ${\cal E}$ into itself, provided $\e$ is chosen small 
enough. Hence, we have obtained a fixed point of the mapping ${\bf C}_\e$. 
This completes the proof of the existence of the hypersurfaces $S_\e$ which 
are described in the  Theorem 1. Most of the properties states in Theorem 1
follow readilly from the construction itself except the derivation of (\ref{eq:slope}). 

\medskip

{\bf Proof of (\ref{eq:slope})~:} This follows from the application of the well known 
balancing formula for minimal hypersurfaces. In the case where $m=n-1$ we know 
from the construction itself that the hypersurface $S_\e$ is, away from the origine, 
the graph of the function 
\[
(x_1, x_2 ) \longrightarrow c_\e \, |x_1| + c_\infty \, \e + 
{\cal O}(\e^{n-1} \, |x_1|^\nu)
\]
Moreover, near $0$ the hypersurface is a graph over the rescalled $n$-catenoid.
It remains to identify the constant $c_\e$. In order to do so, 
we apply the balancing formula of \cite{Ros} (Theorem 7.2) between 
the hyperplane $z = 0$ and $z=z_0$ for $z_0$ tending to $+\infty$. This yields
\[
\mbox{Vol} (T^{n-1} ) c_\e \sim \e^{n-1} \, \mbox{Vol} (S^{n-1}) .
\]
And (\ref{eq:slope}) follows at once from our normalization of
the volume on an $n-1$-dimentional torus.

\end{document}